\documentclass[12pt]{article}
\usepackage{amssymb}
\usepackage{amsmath}
\newtheorem{theo}{Theorem}
\newtheorem{var}[theo]{Variant}
\newtheorem{prop}[theo]{Proposition}

\newtheorem{coro}[theo]{Corollary}
\newtheorem{rema}[theo]{Remark}
\newtheorem{Defi}[theo]{Definition}

\newtheorem{conj}[theo]{Conjecture}
\newtheorem{quest}[theo]{Question}
\newtheorem{Defi-prop}[theo]{Definition-Proposition}
\newtheorem{claim}[theo]{Claim}

\def \F{{\mathbb F}}

\def \N{{\mathbb N}}
\def \A{{\mathbb A}}

\def \R{{\mathbb R}}
\def \Z{{\mathbb Z}}
\def \P{{\mathbb P}}
\def \Q{{\mathbb Q}}
\def \p{{\mathfrak p}}

\begin{document}

\title{Some applications of $p$-adic uniformization to algebraic dynamics}

\author{Ekaterina Amerik}

\date{}

\maketitle

The purpose of these notes is to give a brief survey of several topics
at the limit of geometry and arithmetics, where some fairly elementary 
$p$-adic methods have led to highly non-trivial results. These results
are recent but not brand-new: all proofs have been published elsewhere.
My hope and reason for putting them together is that this might facilitate
further progress, in particular by young mathematicians or those who are new 
to the field. Since I am a geometer which only meets arithmetics by accident,
the point of view in this notes is quite biased. The reader is encouraged to
consult the texts by other people who have contributed to the subject:
in particular, the forthcoming book by Bell, Ghioca and Tucker \cite{bgt-book} 
promises 
to be very interesting.

The notes are written for the proceedings volume of CRM Montreal thematic program 
``Rational points, rational curves and entire holomorphic curves on algebraic varieties'' in June 2013.
During the writing of the notes, I  was also preparing a mini-course on the subject 
for ANR BirPol and Fondation Del Duca meeting ``Groupes de transformations'' in Rennes in June 2014.
I am grateful to the organizers of both activities for giving me this opportunity to speak.
Thanks also to Dragos Ghioca for sending me a preliminary version of \cite{bgt-book} and answering a few questions.

\section{A motivation: potential density}\label{intro-pd}

Let $X$ be a projective variety over a field $K$. 

\begin{Defi} Rational points of $X$ are
{\it potentially dense} over $K$ (or, as one also sometimes says, $X$ is potentially dense over $K$) if there is a finite extension $L$ of $K$
such that the $L$-points are Zariski dense in $X$.
\end{Defi}

The reason for looking at the potential density rather than at the density of $K$-rational
points is that the potential density behaves much better from the geometric
point of view. Indeed, even a plane conic over the rationals can have a
dense set of rational points or no rational points at all; whereas if we look
at the potential density, we may at least hope that the varieties which share
similar geometric properties should be potentially dense (or not potentially
dense) all at once.

If $X$ is {\it rational} (that is, birational to $\P^n$), or, more generally,
{\it unirational} (that is, dominated by $\P^n$) over $\bar K$, then rational
points are obviously potentially dense on $X$. Indeed, choose $L$ such
that the unirationality map $f: \P^n\dasharrow X$ is defined over $L$.
Then the images of $L$-points of $\P^n$ are $L$-points of $X$ and they are
Zariski dense in $X$ since $f$ is dominant and $L$-points are dense on $\P^n$.
More generally, a variety dominated by a potentially dense variety is itself
potentially dense.

It is certainly not true in general that a variety which dominates a 
potentially dense
variety is itself potentially dense. However this is true in an important
particular case: if $f:X\to Y$ is a finite \'etale morphism and $Y$ is
potentially dense, then so is $X$. This follows from Chevalley-Weil theorem
(see for example \cite{serre}). The idea is that points in the inverse image
of $x\in X(K)$ are points over finite extensions of degree equal to $deg(f)$
and ramified only at a fixed (that is, independent of $x$) finite set of places.There is only a finite number of such extensions.

Unirational varieties share many other properties of the projective space. For
instance, the tensor powers of the canonical line bundle $K_X$ on such a variety
$X$ have no sections: indeed, a section of $K_X^{\otimes m}$, $m>0$, 
would pull back to $\P^n$ and give
a section of a tensor power of $K_{\P^n}$ (by Hartogs' extension theorem), 
but no such section exists.

On the opposite geometric end, we have the {\it varieties of general type}: these are
the varieties on which the tensor powers $K_X^{\otimes m}$ have ``lots of sections''.
More precisely, a smooth projective variety $X$ is said to be {\it of general type}
if the map defined by the linear system $|K_X^{\otimes m}|$ is birational to its 
image for some $m>0$. The simplest examples are curves of genus $g\geq 2$, or
smooth hypersurfaces of degree $\geq n+2$ in $\P^n$. The following conjecture is very famous:

\begin{conj}(Lang-Vojta) A smooth projective variety $X$ which is 
of general type cannot be potentially dense over a number field $K$.
\end{conj} 

Up to now this is known only for curves and for subvarieties of abelian varieties,
by the work of G. Faltings.

Lang-Vojta conjecture implies that varieties dominating a variety of general type
cannot be potentially dense over a number field. One might ask whether this should
lead to a geometric characterization of potentially dense varieties. 

The naive guess is wrong: one can construct a surface which is not of 
general type and does not dominate
any curve of genus $g\geq 2$, yet it is not potentially dense, since it
 admits a finite \'etale covering which {\it does} map onto a curve of higher genus,
and the potential density is stable under finite \'etale coverings. This seems to be first observed by Colliot-Th\'el\`ene, Skorobogatov and Swinnerton-Dyer 
in \cite{CSS}. 
The idea is to take an elliptic surface $X$  over $\P^1$ with suitably many double
fibers, so that locally the map to $\P^1$ looks like $(x,y)\mapsto u=x^2$. Then the ramified covering $C$ of the base which eliminates these multiple
fibers (that is, locally looks like $z\mapsto u=z^2$, so that the fibered product is singular
and its normalization is \'etale over $X$)  will have genus at least two. 

F. Campana suggests in \cite{C} that the potentially dense varieties are exactly
the so-called {\it special varieties}. Roughly speaking, these are the varieties
which do not dominate
{\it orbifolds of general type}: if $f:X\to B$ is a fibration with certain
good properties (which are achieved on suitable birational models), one
can define an orbifold canonical bundle $K_B+\Delta$ by taking into account the
multiple fibers of $f$, and this bundle should not have too many section. For the
moment, proving this, or even the ``easier'' direction that special varieties
should be potentially dense, looks quite out of reach. 

In any case, all existing philosophy seems to imply that the varieties with 
negative canonical bundle (the {\it Fano varieties}) or trivial canonical bundle
must be potentially dense. There is a reasonable amount of evidence for this in 
the Fano case: indeed many Fano varieties are known to be unirational, and when
the unirationality is unknown the potential density still 
can sometimes be proved
(see for example \cite{HarT}). Also, potential density is known for tori (and it shall be 
explained in this survey in a particularly elementary way).
But the case of simply-connected varieties with trivial canonical class remains
mysterious: indeed, even for a general K3 surface the answer is unknown, and
moreover there is no example of a potentially dense K3 surface with cyclic Picard
group (that is, ``general'' in the moduli of polarized K3). 

Bogomolov and Tschinkel \cite{BT} proved the potential density of elliptic K3 
surfaces.
\begin{theo}
Let $X$ be a K3 surface over a number field $K$. If $X$ admits an elliptic
fibration, then $X$ is potentially dense.
\end{theo}

{\it Idea of proof:} Construct a multisection $C$ which is a rational curve
(so has a lot of rational points) and
is non-torsion, that is, the difference of at least two of its points on a 
general fiber is non-torsion in the jacobian of this fiber.
Then one can move $C$ along the fibers by ``fiberwise multiplying it by an integer''
and produce many new rational points in this way.

\medskip

More generally, let $X$ be a variety equipped with a rational self-map
$f:X\dasharrow X$, both defined over a number field $K$ (such as the fiberwise 
multiplication by $k$ on an elliptic 
surface;
this exists for any $k\in \Z$ when the surface has a section and for suitable $k$ if not). 
It is a natural idea to use $f$ to produce many rational points on $X$: indeed 
$f$ sends rational points to rational points. 

This approach has first been worked out by Claire Voisin and myself \cite{AV}
to give the first example of a simply-connected variety with trivial canonical
class which has Picard number one (so is ``general'' in the polarized moduli 
space)
and has potentially dense rational points. Our example is as follows. 

Let $V$ be a cubic in $\P^5$ and $X={\cal F}(V)\subset Gr(1,5)$ be the variety
parameterizing the lines on $V$. A simple computation shows that $X$ is a smooth
simply-connected fourfold with trivial canonical bundle. Moreover it can be seen 
as a higher 
dimensional analogue of a K3 surface: as shown by Beauville and Donagi \cite{BD},
$X$ is an irreducible holomorphic symplectic manifold (that is, $H^{2,0}(X)$ is
generated by a single nowhere degenerate form $\sigma$), deformation equivalent
to the second punctual Hilbert scheme $Hilb^2(S)$, where $S$ is a K3 surface
(and actually isomorphic to $Hilb^2(S)$ when the cubic $V$ is pfaffian).

\begin{prop}(C. Voisin) $X$ admits a dominant rational self-map $f:X\dasharrow X$
of degree 16.
\end{prop}

{\it Sketch of proof:} Let us describe the construction: for a general 
line $l$ on $V$,
there is a unique plane $P$ tangent to $V$ along $l$ (indeed the normal 
bundle $N_{l,V}={\cal O}_l\oplus {\cal O}_l\oplus {\cal O}_l(1)$, which makes this
unicity appear on the infinitesimal level). One defines $f(l)$ as the only line
which is residual to $l$ in the intersection $P\cap V$, and one shows (using,
for example, Mumford's trick on algebraic cycles and differential forms) that
$f$ multiplies $\sigma$ by $-2$.

\begin{theo}(\cite{AV}) For ``most'' cubic $4$-folds $V$ defined over a number field, the corresponding
variety $X={\cal F}(V)$ (which is defined over the same number field) has cyclic
geometric Picard group and is potentially dense.
\end{theo}

What is meant by ``most'' can be made precise, but this is a rather complicated
condition. Since it is not related to our main subject, let us only
mention that the parameter point of the cubic fourfold in question should be
outside of a certain thin subset, like in Hilbert
irreducibility theorem.

The proof, too, is long and involved; in fact
most of my contribution to the main subject of these notes grew out of a search for a 
a more elementary argument. Let us only mention the starting point:
we consider a family of birationally abelian surfaces $\Sigma_t, t\in T$
covering $X$ (the existence of such a family was observed by Claire Voisin
in relation to Kobayashi pseudometric issues) and 
remark that since rational points are potentially dense on $\Sigma_t$ for algebraic
$t$, it is enough to find an
algebraic $t$ such that the 
iterates $f^{k}(\Sigma_t)$ are Zariski dense in $X$.

It turns out to be surprisingly difficult to show by the methods of complex geometry
 that the iterates of something
algebraic are Zariski dense. Let me illustrate this point by explaining the 
difference with the transcendental situation. 

The following theorem has been proved by Campana and myself in the complex geometry
setting.

\begin{theo}\label{fibration-ac}
Let $X$ be a projective variety and $f:X\dasharrow X$
a dominant rational self-map, both defined over an algebraically closed field $K$. 
Then there is a dominant rational map $g:X\dasharrow T$, such that $gf=g$ and for a sufficiently general point $x\in X$, the fiber of $f$ 
through $x$
is the Zariski closure of the iterated orbit $O_f(x)=\{f^k(x),k \in Z\}$. 
\end{theo}

One can always Stein-decompose $g$ to arrive to a map with connected fibers preserved by a
power of $f$.
The theorem thus implies that if no power of $f$ preserve a non-trivial rational fibration
(and this is something which often can be easily established by geometric 
methods, see for example \cite{AC}, theorem 2.1 and corollary 2.2), the orbit of a sufficiently
general point is Zariski dense. If, on the contrary, some power of $f$ does preserve a fibration,
then this is obviously not the case. 

Unfortunately ``sufficiently general'' in the theorem means ``outside a countable union of
proper subvarieties'' (the theorem is proved by looking at the Chow components 
parameterizing $f$-invariant subvarieties and discarding the 
families which do not dominate $X$). That is, when the field $K$ is uncountable, most $x\in X$
are indeed general in this sense; but the theorem does not give any information
when $K$ is countable, since it might happen that no $K$-point is sufficiently
general!

In particular, we still do not know whether there are algebraic points on
the variety of lines of a cubic fourfold which have Zariski-dense iterated
orbit under $f$.  What we do know is that $f$ does not preserve a rational
fibration, and neither do its powers, by \cite{AC}, theorem 2.1; but apriori the iterated orbits of algebraic 
points can have smaller
Zariski-closure than those of general complex points.

One would like to conjecture that in reality it never happens: this is already 
implicit in \cite{AC}.

\begin{conj}\label{orbitclosure} Let $X$ be an algebraic variety with a dominant rational self-map
$f:X\dasharrow X$ defined over a number field $K$.
Consider the map $g:X\dasharrow T$ from theorem \ref{fibration-ac}, and let 
$d$ denote its relative dimension. Then there exists an algebraic point $x\in X(\bar 
\Q)$ such that the Zariski closure of $O_f(x)$ is equal to $d$.
\end{conj}

Some less general versions have been formulated by other authors; for instance,
the following conjecture has been made by Shouwu Zhang. For $X$ a smooth 
projective variety, let us call 
an endomorphism $f:X\to X$ {\it polarized}, if there is an ample line bundle $L$
on $X$ such that $f^*(L)=L^{\otimes q}$ with $q>1$.

\begin{conj} (Zhang) Let $f:X\to X$ be a polarized endomorphism of $X$ defined over a
number field $K$. Then there exists a point  $x\in X(\bar \Q)$ with Zariski-dense
iterated orbit $O_f(x)$.
\end{conj}

Note that a polarized endomorphism cannot preserve a fibration. Indeed, otherwise 
let $F$ be a fiber; one then should have $deg(f|_F)=deg(f)$. But the former is
$q^{dim(F)}$ and the latter $q^{dim(X)}$, a contradiction.

Therefore Zhang's conjecture would follow from conjecture \ref{orbitclosure}.
Indeed, since no power of $f$ preserve fibrations, $T$ is a point and if the
conjecture \ref{orbitclosure} is true, there is an algebraic point with 
Zariski-dense orbit.

In what follows, we shall try to explain some $p$-adic ideas towards the proof of 
this conjecture.

One should mention, though, that there is no hope to prove the potential
density of all special varieties using rational self-maps, as the self-maps 
do not
always exist. For instance, Xi Chen \cite{Ch} proved that a general K3 surface
does not admit a non-trivial rational self-map. Nevertheless, an interesting
example (variety of lines of a cubic fourfold) has been studied in this way, and 
hopefully more shall follow. 

\section{Near a fixed point}\label{abr}

While working on problems of holomorphic dynamics, one is often led to consider
the behaviour of the map in a neighbourhood of a fixed point. In \cite{ABR}, we have
tried to work out some rudiments of a similar approach in algebraic geometry in order
to simplify and render more explicit the proof of potential density of the variety
of lines on a cubic fourfold from \cite{AV}. Somewhat later, we have learned that
similar ideas were exploited by Ghioca and Tucker in order to settle a case of the 
so-called {\it dynamical Mordell-Lang conjecture} to which we shall return in the
next section.

Endomorphisms often have periodic points: for instance, a theorem by Fakhruddin 
\cite{F}
asserts that a polarized endomorphism has a Zariski-dense subset of periodic points.
Replacing $f$ by a power if necessary, we may assume that some periodic point
is actually fixed.

If $X={\cal F}(V)$ is the variety of lines of a cubic in $\P^5$ and $f:X\dasharrow X$
is the rational map which sends a general line $l$ to the line $l'$ which is
residual to $l$ in the intersection of $V$ with the plane tangent to $V$ along $l$,
then the fixed points are, obviously, the lines such that there is a plane tritangent
to $V$ along this line (and not contained in the indeterminacy locus, that is, 
this tritangent plane should be the only plane tangent to $V$ along $l$). An 
explicit computation shows that such lines form a surface on $X$, and no component
of this surface is contained in the indeterminacy locus (one
can, for instance, remark that the fixed surface is certainly lagrangian,
because of the identity $f^*{\sigma}=-2\sigma$, where $\sigma$ is the 
symplectic form, and that the indeterminacy locus is not lagrangian because
of the computations in \cite{A0}; but there is probably a much easier way).

\subsection{Linearization in a $p$-adic neighbourhood}

Let $X$ be arbitrary, and let $q$ be a fixed point of a rational map 
$f:X\dasharrow X$. Assume that everything is defined over a number field $K$.
We shall denote by ${\cal O}_K$ the ring of integers, by $\mathfrak p\subset 
{\cal O}_K$
an ideal, by $O_{\mathfrak p}$ and $K_{\mathfrak p}$ the $\mathfrak p$-adic 
completions.
  
Our starting observation is that for a suitable $\mathfrak p$, one can find a 
$\mathfrak p$-adic neighbourhood $O_{{\mathfrak p}, q}$ (that is, the set of
$\mathfrak p$-adic points reducing to the same point as $q$ modulo $\mathfrak p$
in a suitable model of $X$) which is invariant under $f$, and $f$ is well-defined
there. 

One can define and describe 
$O_{{\mathfrak p}, q}$ in a very down-to-earth way, without
talking about models, by a ${\mathfrak p}$-adic version of the implicit function
theorem.

Namely, following \cite{ABR}, choose an affine neighbourhood $U\subset X$ of $q$,
such that the restriction of $f$ to $U$ is regular. By Noether
normalisation lemma, there is a finite $K$-morphism 
$\pi =(x_1,\dots,x_n):U\longrightarrow\A^n_K$ to the affine space, 
which is \'etale at $q$ and which maps $q$, say, 
to 0. Then the $K$-algebra ${\mathcal O}(U)$ is integral over 
$K[x_1,\dots,x_n]$, 
i.e., it is generated over $K[x_1,\dots,x_n]$ by some regular functions 
$x_{n+1},\dots,x_m$ integral over $K[x_1,\dots,x_n]$. We can view 
 $x_{n+1},\dots,x_m$ and $f^{\ast}x_1,\dots,f^{\ast}x_m$ as power series 
in $x_1, \dots, x_n$ with coefficients in $K$ (indeed the coordinate ring of $U$ is embedded into the local ring of $q$ and 
the latter is embedded into its completion). Since everything is 
algebraic over $K(x_1,\dots,x_n)$, one can show that all coefficients lie in a
finitely generated $\Z$-algebra (this goes back to Eisenstein for $n=1$,
see \cite{ABR}, lemma 2.1). In particular, for almost all primes $\p$,
the coefficients of our power series are $\p$-integral. Take such a $\p$ 
satisfying the following extra condition: for $n<i\leq m$, 
let $P_i$ be the minimal polynomial of $x_i$ over $x_1,\dots,x_n$. We want 
$x_i(q)$ to be a simple root of $P_i(q)$ modulo $\p$ (this condition is 
obviously expressed in terms of the non-vanishing of derivatives modulo
$\p$, and thus also holds for almost all $\p$).

Set $$O_{{\p}, q, s}=\{t\in U(K_{\p})|x_i(t)\equiv x_i(q)\pmod {\p^s}\ for
\ 1\leq i\leq m\},$$ 
and let $O_{\p,q}=O_{\p,q,1}$.

View all our functions $x_i$, $f^*x_i$ as elements of 
${\cal O}_{\p}[[x_1,\dots, x_n]]$.
The following properties are then obvious by construction.

\begin{prop} \label{nbhd} (\cite{ABR}, Prop. 2.2)  (1) The functions $x_1, \dots x_n$ give a bijection
between $O_{{\p}, q, s}$ and the n-th cartesian power of $\p^s$.

(2) The set $O_{\p,q}$ contains no indeterminacy points of $f$.

(3) $f(O_{{\p}, q, s})\subset O_{{\p}, q, s}$ for $s\geq 1$. Moreover, 
$f:O_{\p,q,s}\stackrel{\sim}{\longrightarrow}O_{\p,q,s}$ is bijective 
if $\det Df_q$ is invertible in ${\mathcal O}_{\p}$. 

(4) The $\bar \Q$-points are dense in $O_{{\p}, q, s}$.
\end{prop}

Indeed, the surjectivity just means that for any values of 
$x_i, 1\leq i\leq n$ in $\p^s$, the power series $x_j, j>n$ are going to
converge and have the same value as $x_j(q)$ modulo $\p^s$. But this is clear
since their coefficients are in ${\mathcal O}_{\p}$ (recall that a $\p$-adic
power series converges if and only if the $\p$-adic order of its general 
term goes to infinity). The injectivity is assured by the simple root condition
above: the minimal polynomial $P_i$ evaluated at a point 
$(x_1^0,\dots,x_n^0)\in (\p^s)^{\times n}$
has several roots, but only one of them is congruent to $x_i(q)$ modulo $\p$,
so that only one point out of all pre-images of $(x_1^0,\dots,x_n^0)$ by $\pi$
 is in the corresponding $\p$-adic neighbourhood.

 The power series determining $f$ have integral coefficients and their
constant term vanishes (since $0$ is a fixed point); the second and third
assertions follow easily. Finally the density of $\bar \Q$-points follows from 
their density in $\p^s$.

\medskip

Assume from now on that our fixed point is {\it non-degenerate}, that  is, the
differential $D_qf$ is invertible.
Our next observation is that under a certain condition on the eigenvalues
of the differential at $q$, the map $f$ admits a {\it linearization} in a
small neighbourhood of $q$, complex as well as $\p$-adic. This linearization
is of course analytic, not algebraic, and is a standard tool in holomorphic
dynamics. 

It is usually done in two steps:

\smallskip

\noindent {\bf Step 1} Write a formal power series $h$ conjugating $f$ to its 
linear part $\Lambda$; that is, such that $h\Lambda=fh$.

\smallskip

\noindent {\bf Step 2} Prove that this power series has non-zero convergence radius.

\smallskip

The first step is classical;
one can find it in the textbooks like \cite{arnold}. Let $\lambda_1, \dots,
\lambda_n$ be the eigenvalues of $D_qf$. When
one starts formally solving the equation $f(h(x))=h(\Lambda(x))$, writing
down $h$ term-by-term, one gets the expressions of the form 
$\lambda_1^{m_1}...\lambda_n^{m_n}-\lambda_j$ in the denominators. Here 
$m_i\in \Z,\ m_i\geq 0$ and $\sum m_i\geq 2$. The reason for the last
unequality is that $m=\sum m_i$ is the degree of the corresponding term of $h$ 
(and the term of degree one may be taken to be identity so there is nothing 
to solve).

\begin{Defi-prop}\label{resonance} The equality of the form 
$\lambda_1^{m_1}...\lambda_n^{m_n}-\lambda_j=0$
with $m_i$ as above is called a {\it resonance}. In the absence of resonances,
the map $f$ is formally linearized at $q$.
\end{Defi-prop} 
 
\begin{rema} Note that since we are concerned only with the case when $\sum m_i\geq 2$,
the equality $\lambda_i=\lambda_j$ is not a resonance.
\end{rema}

As for the second step, it turns out to be automatic in the case when 
everything is defined over a number field.

\begin{theo}\label{linearization}(Siegel, Herman-Yoccoz \cite{HY}, Baker, Yu \cite{Yu})
The formal power series $h$ linearizing $f$ has non-zero radius of
convergence as soon as 
$$|\lambda_1^{m_1}...\lambda_n^{m_n}-\lambda_j|>Cm^{-\alpha}$$
for some positive constants $C, \alpha$, and $m=\sum m_i$.
The norm is the usual norm in the complex setting and the
$p$-adic norm in the $\p$-adic setting. This 
condition is always satisfied when the $\lambda_i$ are non-resonant algebraic
numbers.
\end{theo}

This is the classical Siegel condition in the complex setting. The $p$-adic
version is due to Herman and Yoccoz. It follows from the work of Baker that
algebraic numbers satisfy Siegel condition for the complex norm, and the
$p$-adic estimate follows from the work of Yu.

\smallskip

Using the $\p$-adic linearization, we obtain the following sufficient condition
for the density of iterated orbits:

\begin{prop}\label{indep-dense} Suppose that the eigenvalues  $\lambda_1, \dots, \lambda_n$ of $D_qf$ are multiplicatively
independent. Then any point in a suitable $\p$-adic neighbourhood, lying outside of the union of the
coordinate hyperplanes in the linearizing coordinates, has Zariski-dense iterated orbit. In particular,
there is such a point in $X(\bar \Q)$. Therefore rational points are potentially dense on $X$.
\end{prop}

The proof proceeds by linearizing the map and showing that the coordinates of the 
iterates $(\lambda_1^kx_1, \dots, \lambda_n^kx_n), k\in \N$ do not satisfy any analytic equation simultaneously for all $k$
(whereas if they were not Zariski dense, a polynomial function on $X$ vanishing at all iterates
would give a convergent power series in the linearized coordinates).
This relies on an independence property of $\p$-adic exponentials (see lemma 2.6 of \cite{ABR}),
which, as far as I know, is specific for $p$-adics.

Notice that the multiplicative independence, and even the absence of resonances, is only possible
when the fixed point $q$ is isolated. When it is not, then there is a whole subvariety of fixed points,
and the eigenvalue of the differential along this subvariety is equal to one. This means that
the eigenvalues are {\it always} resonant.

However, in the case when ``all resonances are explained by this'', one can still prove a 
linearization statement.

Namely, take $q$ to be a smooth point of the fixed point locus, and let $r$ denote the dimension 
of the fixed subvariety at $q$. One may suppose that $F$ is given by $\{x_{r+1}=\dots=x_n=0\}$.
The first $r$ eigenvalues are then equal to 1. We assume that there are no resonances among the 
other
eigenvalues. In this situation, we have obtained in \cite{ABR} linearization statements
under some technical conditions, for example:

\begin{prop}\label{linear-still}  

In the situation as above, 
suppose that the tangent map $D_qf$ is semisimple and that its
eigenvalues of $D_qf$ do not vary with $q$. 
Under the non-resonance condition
$\lambda_{r+1}^{m_{r+1}}\cdots\lambda_n^{m_n}\neq\lambda_i$ for 
all integer $m_{r+1},\dots,m_n\ge 0$ with $m_{r+1}+\dots+m_n\ge 2$ 
and all $i$, $r<i\le n$, 
the map $f$ can be linearized in some 
$\p$-adic neighbourhood $O_{\p,q,s}$ of $q$, i.e., there exists a formally 
invertible $n$-tuple of formal power series $h=(h^{(1)},\dots,h^{(n)})$ 
in $n$ variables $(x_1,\dots,x_n)=x$ convergent together 
with its formal inverse on a neighbourhood of zero such that 
$h(\lambda_1x_1,\dots,\lambda_nx_n)=f(h(x_1,\dots,x_n))$. \end{prop}

The analogue of the proposition \ref{indep-dense} in this situation is that many algebraic points
have orbit closure of dimension at least equal to the number of the multiplicatively independent
eigenvalues (that is, to the rank of the multiplicative group generated by the $\lambda_i$).
Notice that this is again the statement about the $\p$-adic analytic closure: since the
linearization map is not algebraic, it is difficult to make conclusions about Zariski closure
while working in the linearizing coordinates. In most cases, the Zariski closure should be greater
than that (analytic subvarieties are seldom algebraic). But it is not clear how to find a general
method handling this kind of issues: certainly a challenging task. A couple of examples is considered
in the next section.

\subsection{Examples}

As a first example, let us prove that rational points on an abelian variety $A$
are potentially dense (this is well-known and there are other proofs of this
fact, see e.g. \cite{HT}).

\medskip

{\bf Example 1} Any abelian variety is isogeneous to a product of simple abelian varieties,
and potential density is stable under finite \'etale coverings (by Chevalley-Weil theorem; recall the general discussion of potential density in the first section).
So it is enough to prove the statement for a product of simple abelian varieties
$A=A_1\times\dots\times A_n$.
Consider the map $f=f_1\times\dots\times f_n$, where $f_i$ is the multiplication by $k_i$
on $A_i$ and the $k_i$ are multiplicatively independent. We claim that there is a point
in $A(\bar \Q)$ with Zariski-dense iterated orbit. Indeed, linearize $f_i$ in a suitable
$\p$-adic neighbourhood $O_i$ of zero (where $\p$ is the same for all $i$): this is possible
since there are no resonances. Work in the linearising coordinates and take a point
$(x_1, \dots, x_n)$ with nonzero $x_i$. The analytic closure of the set of iterates
is $l_1\times\dots\times l_i$ where each $l_i$ is a line generated by $x_i$. Now induction
by $n$ and simplicity of the $A_i$ show that $l_1\times\dots\times l_i$ is not contained in a
proper abelian subvariety. But it is easy to see that the Zariski closure
of each $l_i$ must be an abelian subvariety (because of the invariance by
multiplication), 
so the Zariski closure of $l_1\times\dots\times l_i$ must coincide with $A$.

\medskip

{\bf Example 2} Let $X={\cal F}(V)$ be the variety of lines on a cubic fourfold
and $f: X\dasharrow X$ the rational map considered above.
Recall that it multiplies the symplectic form by $-2$ (and is therefore of degree 16).
We have seen that it has a surface $F$ of fixed points. From $f^*\sigma =-2\sigma$, it
is easy to see that $D_qf$ is semisimple and to compute its eigenvalues at a general point $q\in F$: 
these are $(1,1,-2,-2)$. So that we can linearize the map near $q$, but apriori we shall only be
able to conclude that the orbits of points have at least one-dimensional Zariski closure. In \cite{ABR},
we have used the basic idea from \cite{AV} which, together with the linearization
arguments, allowed us to reprove the potential density in a much more
elementary way than in \cite{AV}. Namely, as we have already said, $X$ is covered by a family of birationally abelian 
surfaces: to be precise, these are the surfaces formed by the lines contained in a hyperplane section
of $X$ with three double points. It suffices to prove that the closure of the union
of iterates of such a surface is Zariski-dense (since those which are defined over a number field are
potentially dense). One first shows that the number of such iterates is infinite, by using the $\p$-adic linearisation: in fact what one proves is a general fact that the Zariski 
closure of the 
union 
of iterates of an irreducible subvariety $Y$ passing near a fixed point $q$ remains irreducible as soon
the eigenvalues generate a torsion-free subgroup of ${\mathcal O}_{\p}^{\times}$; 
so that if our surface is not invariant by $f$ (and in our case this is easily checked by
geometry), its iterates are dense at least in a divisor. 
Then one must exclude
the case when the closure of the union of iterates is a divisor $D$: this is again done by 
geometric methods, using case-by-case analysis
on the Kodaira dimension of $D$.

This last step is done in \cite{ABR} assuming that $Pic(X)=\Z$. This holds for most $X$ defined over a 
number field, and there is a technology available to produce explicit examples. The reason is that
$H^2(X,\Q)$ is isomorphic as a Hodge structure to $H^4(V,\Q)$. The space of Hodge cycles in $H^4(V,\Q)$
is known to be one-dimensional for a general complex $V$ (by a version of the ``Noether-Lefschetz theorem''),
but also for most cubics defined over a number field and even over the rationals, namely, for cubics whose 
parameter point lies
outside of a certain thin subset in the parameter space $\P^N$ (\cite{T}; more precisely, this concerns
one-dimensionality of the space of algebraic cycles, but the Hodge conjecture is known for cubic
fourfolds, so this applies to the Hodge cycles as well). It should be possible to get rid of this assumption
by a more careful analysis of the geometry related to the map $f$.

\section{Uniformization of orbits and applications}

\subsection{Dynamical Mordell-Lang problem}

The subject of this subsection is, mainly, the recent work of Bell, Ghioca
and Tucker which has initiated an approach to algebraic dynamics based on the
$p$-adic uniformization of orbits. In their book \cite{bgt-book}, these issues
are going to be explained in much more detail; I shall give below a 
fairly brief presentation.

The following conjecture is known as the dynamical Mordell-Lang problem.

\begin{conj} Let $X$ be a quasiprojective complex variety and $f:X\to X$
an endomorphism. Let $x\in X$ be a point and $V\subset X$ a subvariety.
Then the set of $n\in \N$ for which $f^n(x)\in V$ is a union of a finite 
set and finitely many arithmetic progressions.
\end{conj}

To explain the terminology, recall that the Mordell-Lang conjecture (now a
theorem) reads as 
follows.

\begin{theo} (Faltings for finitely generated groups, McQuillan in general) 
Let $A$ be an abelian 
 variety, $V\subset A$ a subvariety,
$\Gamma\subset A$ a finite rank subgroup. Then $\Gamma\cap V$ is a
finite union of cosets of subgroups of $\Gamma$.
\end{theo}

When $\Gamma$ is the group of torsion points of $A$, this implies the Manin-Mumford
conjecture first established by Raynaud. One may also take for $\Gamma$
the group of rational points $A(K)$, where $K$ is a number field, as 
Mordell-Weil
theorem states that it is finitely generated, and deduce Lang's
conjecture for subvarieties of abelian varieties. Both deductions require 
some work and are based on the observation that an irreducible 
subvariety of an abelian variety which is invariant by multiplication by 
integers must be a subtorus; as a consequence, irreducible components of 
the Zariski closure of a subgroup of $\Gamma$ are subtori, too (they can
of course eventually be trivial, i.e. reduced to points).

 The case when $\Gamma$ is infinite
cyclic is an old result of Chabauty: $\Gamma\cap V$ is a union of a finite set 
and finitely many
arithmetic progressions in $\Gamma$. Sometimes one regards a point in the
finite set as a ``trivial arithmetic progression'' and in this way 
eliminates the finite set
from the statements.

The word ``dynamical'' in this context means that one would like to 
extend this to the 
setting when $X$ is a variety with, say, a finitely generated monoid of 
commuting
endomorphisms, replacing the finitely generated subgroup $\Gamma$ by the
set of points of the form $f_1^{n_1}\dots f_r^{n_r}(x)$ where $f_i$ are the
generators, and trying to understand
the structure of the set of $(n_1,\dots, n_r)$ for which 
$f_1^{n_1}\dots f_r^{n_r}(x)\in V$. This does not look realistic for many
reasons, but the case when $r=1$ (``Chabauty's case'') is especially promising as it reminds one of a 
classical theorem of Skolem-Mahler-Lech about zeroes in recurrent sequences.

Let us recall this story very briefly.

\begin{Defi} A linear recurrence is a sequence of complex numbers 
$(a_n)_{n\in \N}$
such that $a_{n+r}=c_0a_n+c_1a_{n+1}+\dots+c_{r-1}a_{n+r-1}$ for a certain $r$,
constants $c_0,\dots, c_{r-1}$ and any $n$. 
\end{Defi}

\begin{quest} What can one say about the set of $n$ such that $a_n=0$?
\end{quest}

It is well-known that there is an explicit formula for $a_n$:
$a_n=\sum_{i=1}^{k}f_i(n)z_i^n$, where $z_i$ are the roots of the characteristic polynomial and $f_i$ are
polynomials which account for multiple roots; in particular, if the characteristic polynomial does not have
multiple roots, then $r=k$ and $f_i$ are constant. In order to use this formula to solve the problem,
embed everything into $\Q_p$ for a suitable $p$, in such a way that the coefficients of $f_i$ are $p$-adic integers
and $z_i$ are $p$-adic units (by a lemma of Lech \cite{L}, there are infinitely many of such $p$).
The functions $z_i^n$ are not necessarily analytic (that is, given by convergent power series) on $\Z_p$ 
as functions of $n$,
since $\log z_i$ only makes sense when $z_i\equiv 1 \pmod p$. However, $z_i$ is a $p$-adic unit and $p$-adic units
form a finite group modulo $p$. Hence there exists a natural number $N>0$ such that $z_i^N\equiv 1\pmod p$ for
all $i$. The functions $z_i^{Nn}$, and therefore $a_{Nn}$, are analytic in $n$ on $\Z_p$ 
(the series for $\exp$ converges when the $p$-adic order of the argument is greater than $\frac{1}{p-1}$, and this is always true for $\log (1+pa)$ with 
integral $a$). The same 
applies to $a_{Nn+m}$ for $0\leq m\leq N-1$.

Now if  $a_{Nn+m}$ vanishes for infinitely many $n$, it must vanish identically, being analytic on a compact $\Z_p$.
Therefore we obtain

\begin{theo}(Skolem, Mahler, Lech) Let $(a_n)_{n\in \N}$ be a linear recurrence. Then the set of $n$ such that
$a_n=0$ is a union of a finite set and a finite number of arithmetic progressions.
\end{theo}

The work by Bell, Ghioca and Tucker puts this into the geometric context to prove the dynamical Mordell-Lang
conjecture for \'etale endomorphisms.

\begin{theo}\label{etdml}(\cite{BGT})
Let $X$ be a quasiprojective variety with an \'etale endomorphism $f: X\to X$, $x\in X$, and $V\subset X$ is a subvariety. Suppose that the iterated orbit
$O_f(x)$ hits $V$ infinitely often. Then there exist $l$ and $k$ such that 
$f^{l+nk}(x)\in V$ for any $n\in \N$.
\end{theo}

As easily seen, this is just a reformulation of dynamical Mordell-Lang in the \'etale case. The 
main idea
of proof is as follows.
We may assume, by the above-mentioned embedding lemma of Lech, that all parameters are $p$-adic integers.
Assume for the moment that we can prove the following claim.

\begin{claim}\label{uniformiz}There exists a number $N$ and $p$-adic analytic 
maps $\theta_i$ from $\Z_p$ to $X$,  such that 
$f^{Nk+i}(x)=\theta_i(k)$.
\end{claim}

The theorem then follows easily. Indeed, suppose that the iterated orbit of $f$ hits $V$ infinitely
often, and let $G$ be one of the equations defining $V$ in $X$. Then for some $i$, the analytic function 
$G\theta_i$ has
infinitely many zeros and therefore must be zero identically; so the set of $n$ for which $f^n(x)\in V$ 
contains an arithmetic progression.

\medskip

J. Bell (\cite{B}) has first constructed such $\theta_i$ in the case when $X$ is affine and $f$ is an automorphism,
using special properties of polynomial automorphisms of the affine space (like constant Jacobian determinant).
Later, Bell, Ghioca and Tucker have realized that this is possible for any \'etale endomorphism of a quasiprojective
variety. The crucial point in their argument was later greatly simplified by B. Poonen (\cite{P}). Let me try
to explain their argument now.

\medskip

\noindent {\bf Step 0}
The starting very basic remark is that without loss of generality we can replace $f$ by a 
power and our 
point $x$ by some iterate $f^k(x)$.

\smallskip

\noindent {\bf Step 1}
The second remark is that for a suitable $p$ (and there are in fact infinitely
many of such), we can embed the ring generated by all our parameters (that is, the coefficients
of the equations defining $X$, $f$ and $V$ and the coordinates of $x$) in $\Z_p$ in such a way as to 
obtain a $\Z_p$-model of $X$ with good properties (smooth and projective, equipped with a $\Z_p$-point,
a subvariety and an unramified map $\Phi$ which models $f$). This is done by some standard yoga 
in algebraic
geometry for which we refer to the original paper \cite{BGT}; the reader not comfortable with these 
notions may use the down-to-earth description of the last section with $\Q_p$ and $\Z_p$ instead of
$K_{\p}$ and ${\cal O}_{\p}$.

\smallskip

\noindent {\bf Step 2}
The next observation is that, replacing if necessary the map $f$ by a power,
and $x$ by an iterate, we can construct an invariant $p$-adic 
neighbourhood $O_{p, x}$ of $x$ as in the last section (more precisely,
we have in mind $O_{p, x, 1}$ from the last section, i.e. points reducing
to the same point as $x$ modulo $p$) and provide a similar description of this neighbourhood.

Indeed, to start with, let us forget about the endomorphism. Then one may 
either try to reason as in the last section (choosing a $p$ in a
right way), or to remark that points of the $\Z_p$-model ${\cal X}$ which
reduce to the same point as $x$ modulo $p$ correspond to the prime ideals
in ${\cal O}_{{\cal X},x}$, and in fact also in its completion which is just
the ring of power series in $n$ variables $x_1,\dots, x_n$, such that the
residue ring is $\Z_p$. The latter are of the form $(x_1-pb_1,\dots, x_n-pb_n)$
(where $b_i\in\Z_p$).

After that, we define $O_{x, q}$ as in the last section and identify it
with $(p\Z_p)^n$ in the same way.

\smallskip

\noindent {\bf Step 3}
The slight difference is where the map $f$ (or $\Phi$) is concerned.

Indeed, recall that in the last section we have constructed an invariant $p$-adic neighbourhood
of a {\it fixed} point $q$. Our point $x$ is not fixed, but it is certainly preperiodic modulo $p$,
since a (quasi)projective variety over a finite field has only finitely many points (and thus all orbits
modulo $p$ are finite). Replacing $x$ by an iterate and $\Phi$ by a power, we may assume that $x$ is
fixed modulo $p$. 

The functions defining $\Phi$ are then power series in $x_1,\dots, x_n$ with integral coefficients and
the constant term divisible by $p$ (recall that it was zero in the last section). As before, the map
is well-defined on  $O_{p, x}$.

It shall be convenient for us to identify  $O_{p, x}$ with $\Z_p^n$ rather than with  $(p\Z_p)^n$.
We have therefore to make a change of variable, writing $H(X)=\frac{1}{p}\Phi(pX)$, where
$X=(x_1,\dots, x_n)$. In this way,
 $\Phi$ is given on $\Z_p^n$ by the power series $H_1,\dots, H_n$ with $\Z_p$-coefficients, converging on
the whole of $\Z_p^n$. Moreover, the terms of degree at least two in $H$ are zero modulo $p$:
$HX\equiv C+LX \pmod p$ where $C$ is constant and $L$ is linear.  

More generally, the coefficients by terms of degree $m$ must be divisible
by $p^{m-1}$.

\smallskip

\noindent {\bf Step 4}
It is time to explain where the unramifiedness condition comes in. We actually only need that $f$ is unramified
modulo $p$. If such is the case, the differential of $f$ is nowhere degenerate; one immediately deduces that
the linear map $L$ is non-degenerate. Therefore $X\mapsto C+LX$ is an 
automorphism of a finite-dimensional $\F_p$-affine space, which is a
finite set. Therefore some power of it is the identity. Hence, again replacing 
our endomorphism by a power,
we may assume that $H$ is identity modulo $p$.

Now the claim (and thus the theorem \ref{etdml}) clearly follows from the proposition below.

\begin{prop}\label{bell} (Bell) 
Let $H:\Z_p^n\to \Z_p^n$ be a $p$-adic (where $p>3$) analytic map given by power series
with integral coefficients, and suppose that $H(X)\equiv X\pmod p$. Assume 
a mild growth condition: the order of the coefficient of a term of degree $m$
is at least $m-1$ (we have just remarked above that it is satisfied in the
situation that we want to understand).  Let 
$x\in \Z_p^n$.
Then there exists a $p$-adic analytic map $G: \Z_p\to \Z_p^n$ (given by
power series with $\Q_p$-coefficients but converging on the whole of $\Z_p$) 
such that
$G(n)=H^n(x)$ for $n\in \N$.
\end{prop}

Poonen gave an especially simple and efficient variant of this in \cite{P}.

\begin{var}\label{poonen} Let $K$ be a field complete with respect to the absolute value
satisfying $|p|=1/p$. Let $R$ be the valuation ring. Denote by $R\langle X\rangle$ the 
Tate algebra (of power series converging on the closed unit polydisc, i.e.
those whose coefficients tend to zero as the degree grows). If $F\in R\langle X\rangle ^n $
satisfies $F(X)=X\pmod {p^c}$ for some $c>\frac{1}{p-1}$, then there exists
$G\in R\langle X,t\rangle ^n$ with $G(X,n)=F^n(X)$ for $n\in\N$. 

\end{var}

{\it Sketch of proof:} Notice that the Tate algebra is the completion of
the polynomial ring with respect to the norm which is the maximal absolute 
value
of the coefficients. Consider the difference linear operator on analytic maps from $R \langle X\rangle ^n $:
$$\Delta H(X)= H(F(X))-H(X).$$ Since $F(X)=X\pmod {p^c}$, one infers that it maps 
$R[X]$ to $p^cR[X]$ and therefore $R \langle X\rangle ^n $ to $p^cR \langle X\rangle ^n $. That is, 
$\Delta^mX\in p^{mc}R\langle X\rangle ^n$. This implies convergence of the Mahler series 
$$G(X, t)=\sum_{m\geq 0} t(t-1)\dots(t-m+1)\frac{\Delta^mX}{m!},$$ 
and it is easy to check that $G(X,n)=F^n(X)$.

\medskip

\begin{rema}\label{numberfd} One useful remark (\cite{P}, remark 4) is that if $F(X)=X\pmod {p^c}$ for some $c>0$,
then $F^p(X)=X\pmod {p^{c'}}$ for a larger $c'$, and in this way we can increase $c$ to be
greater than $\frac{1}{p-1}$ after replacing $F$ by a power. This makes possible
to obtain and use an analogue of claim \ref{uniformiz} and related statements over 
a finite extenstion $K_{\p}$ of $\Q_p$, in an obvious way. This observation has been 
made earlier in \cite{A} (theorem 7), following the original argument of Bell.
\end{rema}

\subsection{Points with infinite orbit}

We are going to apply the $p$-adic uniformization of orbits to a problem of completely different flavour.
Recall the ``generalized Zhang conjecture'' from the section \ref{intro-pd}: it affirms
that if $X$ is a projective variety over a number field $K$ and 
$f: X\dasharrow X$ is a dominant rational self-map, then one can 
always find an algebraic point such that the dimension of its
iterated orbit is as large as possible (that is, coincides with the relative
dimension $d$ of the associated fibration $f: X\dasharrow T$, see 
theorem \ref{fibration-ac}).

Notice that apriori it is not even clear that the iterated orbit is 
well-defined for certain $x\in X(\bar \Q)$. Indeed, apriori it could happen
that for any $x\in X(\bar \Q)$, there is a number $n$ such that $f^n(x)$ is
in the indeterminacy locus of $f$. The union of the inverse images of the
indeterminacy locus is a countable union of proper subvarieties, but $\bar \Q$ itself
is countable, and so is $X(\bar \Q)$.

Fortunately this turns out to be true. But the proof I know appeals to a difficult theory by
Hrushovski (\cite{H}), which implies that any rational self-map has a lot of periodic points
over finite fields.

Let me state Hrushovski's main result from \cite{H}.

\begin{theo} Let $X$ be an affine variety of dimension $d$ over $\F_q$ and 
$S\subset X\times X$ 
a correspondence
such that both its projections to $X$ are dominant and one of them is quasi-finite. Let 
$\Phi_q$ denote the Frobenius map. Then for $m>>0$, the intersection 
$(S\cap \Gamma_{\Phi_{q^m}})(\bar \F_q)$ is non-empty. More precisely, 
$|(S\cap \Gamma_{\Phi_{q^m}})(\bar \F_q)|=
aq^{md}+O(q^{m(d-1/2)})$ for $m>>0$ and a constant $a$.
\end{theo}  

This type of estimate has been known before, and under stronger hypotheses on $S$ 
more precise results are available (Deligne's conjecture, see for example \cite{Var} and the references therein).
The feature of Hrushovski's theory is exactly that the hypotheses on the correspondence $S$ are
extremely faible: indeed the graph of any dominant rational map would do. In the results of this type
which have been known before, one asks, roughly speaking, that one projection is proper and 
another is quasifinite. This is very rarely achieved when $S$ is the graph of a rational, but
not regular endomorphism: as soon as one removes some subvariety from $X$ to render the second 
projection quasifinite, the first one ceases to be proper. Here, on the contrary, one may remove 
any subvariety from $X$ without destroying the conditions, and therefore still have periodic points.

\begin{coro}  Let $X$ be an algebraic variety of dimension $d$ over $\F_q$ and $S\subset X\times X$ 
the graph of a dominant rational self-map. For any proper subvariety $Y\subset X$, there exists an $m>>0$
such that $(S\cap \Gamma_{\Phi_{q^m}})(\bar \F_q)$ contains a point $(x, \Phi_{q^m}(x))$,
where $x$ is outside of $Y$.
\end{coro}

Now it is clear that not every point of $X(\bar \F_p)$ eventually lands in the indeterminacy locus.

\begin{coro} In the situation as above, denote by $f$ the rational self-map of $X$ which has
$S$ as a graph and set $Y=\cup_i f^i(Indet(f))$. Note that this is a finite union since everything
is defined over a finite field. Then there exists a point $x\in \F_{q^m}$ such that 
$f^i(x)\not\in Indet(f)$ for any $i$; moreover $x$ is periodic for $f$.
\end{coro} 

{\it Proof:} Set  $Y=\cup_{i\in Z} f^i(Indet(f))$. Then by the previous corollary, there is a point $x$ outside $Y$
with the property that $f(x)=\Phi_{q^m}(x)$. One deduces that $f^k(x)=(\Phi_{q^m})^k(x)$
(indeed $f$ commutes with Frobenius since the latter is a field automorpism),
clearly $x$ must be periodic.

\medskip

Finally, such points lift to algebraic points of $X$ (not necessarily periodic).

\begin{coro} Let $f:X\dasharrow X$ be a dominant rational self-map of an algebraic variety
defined over a number field. Then there are algebraic points of $X$ with well-defined
orbit (i.e. such that no iterate of these land in the indeterminacy locus). 
\end{coro}

It is clear that
in this statement (as well as in the previous one) one may replace the indeterminacy 
locus by any other proper subvariety,
for instance the ramification, or the union of ramification and indeterminacy.

For suitable $\p$, the $\p$-adic neighbourhoods $O_{\p,x}$ of such points
are invariant by a power of $f$ (since $x$ is fixed by this power of $f$ modulo $\p$).

This is the key observation needed to make the first step towards the conjecture on
orbit closures, which is the main result of this section.

\begin{theo}\label{nonpreper} Let $X$ be a variety and $f:X\dasharrow X$ a dominant rational self-map
of infinite order, both defined over a number field $K$. Then one can find a $\bar \Q$-point
of $X$ with infinite orbit; moreover the set of such points is dense in $X$.
\end{theo}

In other words, if the Zariski closure of a general complex point has strictly positive dimension,
the same is true for most algebraic points. 

\medskip

{\it Sketch of proof:} Take an algebraic point $x$ as above, periodic modulo 
$\p$ for $\p$ a ``good enough'' prime. 
Replacing $f$ by a power, we may assume that $O_{\p,x}$ is invariant by $f$, which is given by
power series with $\p$-integral coefficients on it. Moreover, if $x$ is not a ramification point
modulo $\p$, one has, in coordinates $X=(x_1,\dots, x_n)$ on $O_{\p,x}$, that $f(X)=C+LX\pmod \p$ 
with $L$ linear non-degenerate, so that if we replace $f$
by a further power, we get $f(X)=X\pmod \p$. Proposition \ref{poonen} does not apply directly
since this is a weaker assumption that $f(X)=X\pmod {p^c}$ with $c>\frac{1}{p-1}$. Indeed, $\p$
is generated by an uniformizer $\pi$ of $p$-adic order $\frac{1}{e}$, where $e$ is the 
ramification index which does not have to be bounded by $p-1$. However, after replacing $f$ by
yet another power, we increase the exponent $c$ of $p$ in the congruence $f(X)=X\pmod {p^c}$ (see 
remark \ref{numberfd}). Therefore by proposition \ref{poonen} 
we can uniformize orbits.
In conclusion, for any $y\in O_{\p,x}$ 
there are $l\in \N$ ($l\neq 0$) and $G_y$ anaytic on $\Z_p$, such that $(f^l)^n(y)=G_y(n)$.

Assume now that $y$ is preperiodic. Then $(f^l)^n(y)$ and therefore $G_y(n)$ takes only finitely
many values. By analyticity, $G_y(n)$ must be constant on $\Z_p$. Therefore $(f^l)^n(y)$ is constant
as $n$ varies. In other words, every preperiodic point in $O_{\p,x}$ is periodic of bounded period (at most $l$).
But such points lie on a finite union of proper subvarieties and there are algebraic points
of $O_{\p,x}$ outside of these subvarieties. This finishes the proof.

\begin{rema} When $f:X\to X$ is regular and polarized, it is easy to prove the existence of non-preperiodic
algebraic points using the canonical heights. Recall that for a regular
polarized map $f: X\to X$, $f^*L=L^{\otimes q}$, $L$ ample, $q>1$, one can 
introduce a function $h: X(\bar \Q)\to \R_{\geq 0}$ which is a Weil height
for $L$ and satisfies $h(f(x))=qh(x)$ (\cite{CS}). Clearly preperiodic points
have zero canonical height, but at the same time $X(\bar \Q)$ is not a set of
bounded height, so there must be non-preperiodic algebraic points. 

So our theorem is not a nontrivial step towards the
proof of Zhang conjecture, but rather towards the more general orbit closure conjecture 
\ref{orbitclosure}.
\end{rema}

\subsection{Bounding periods}

When one tries to extend theorem \ref{nonpreper} further, one sees that the problem splits in two parts,
one dynamical and one geometric. Surprisingly, the dynamical part is easy and it is solved in the
same way as in the case of dimension zero (as it is remarked by Bell, Ghioca and Tucker in \cite{BGT2},
theorem 1.1). The geometric part is difficult and for the moment one
does not seem to have any clue about the general case.

Let me state these two problems precisely. Let us say we want to solve Zhang's conjecture for surfaces:
take a surface $S$ with an endomorphism $f$ such that no power of $f$ preserves a fibration, and look for
an algebraic point $x$ with Zariski-dense orbit. By definition, if the orbit is not Zariski-dense, it is 
contained in a proper subvariety (eventually reducible). This subvariety can be chosen to be
invariant by $f$; in particular, its components of maximal dimension are then periodic. We should therefore
aim to find an algebraic point not contained in a periodic subvariety.

The two parts are as follows.

\noindent {\bf Part 1} Prove that the relevant preperiodic subvarieties
are periodic of bounded period; in particular, replacing $f$ by a power
we may assume that they are invariant.

\medskip

\noindent {\bf Part 2} Prove that, unless $f$ preserves a fibration, 
there is not enough invariant subvarieties
to contain all algebraic points. 

\medskip

As we have already told, Part 1 is completely general and very easy.
The proof of the following observation (see e.g. \cite{BGT2}) 
is the same as the corresponding argument in 
theorem \ref{nonpreper}.

\begin{prop} In the invariant $\p$-adic neighbourhood $O_{\p,q}$ 
constructed as above
and such that a power of the map $f$ admits uniformization of orbits, periodic 
subvarieties $Y$ have bounded period.
\end{prop}

We have already seen that the uniformization condition is satisfied 
when $f$ is \'etale, or, more generally, when the orbit of $q$ does not
intersect the ramification modulo $\p$ -- by Hrushovski theorem we can
always find a $q$ like that, eventually after a finite extension our base field $K$.

\medskip

{\it Proof:} Apriori, one has $f^N(Y)=Y$, where $N$ depends on $Y$. 
Choose $x\in Y$ in $O_{p,q}$, then $f^i(x)\in Y$ infinitely many times
(for $i$ a multiple of $N$, for example). But there is an analytic function
$G_x$ such that $G_x(n)=f^{ln}(x)$. For any polynomial $P$ determining $Y$
in $X$, the $\p$-adic analytic function $PG_x$ has infinitely many zeroes
and therefore is zero identically. We conclude that $f^l(x)\in Y$, and
therefore, taking all possible $x\in Y$ and considering Zariski closure, 
that $f^l(Y)\subset Y$. The proposition is proved.

\medskip 

Part 2 is difficult for the following reason. Whereas points on an algebraic
variety are ``parameterized'' by something irreducible - namely the variety
itself, the higher-dimensional subvarieties are parameterized by the Chow
variety which has countably many components (e.g. the degree of a subvariety
can grow indefinitely). So one can, apriori, have an unbounded family of 
invariant subvarieties in the neighbourhood of a point. And this is what
sometimes indeed happens, e.g. for the $d$th power map 
$(x:y:z)\mapsto (x^d:y^d:z^d)$
on $\P^2$: through the point $(1:1:1)$, there is an infinite number of
invariant curves of growing degree of the type $x=y^{\alpha}$, 
$\alpha\in \Q$ (the degree is $max (|c|, |d|)$ where $\alpha=c/d$ and
$c,d$ are relatively prime integers).

\medskip

For birational maps, there are some positive results. For instance, Serge Cantat
proved in \cite{Can} that the number of invariant hypersurfaces of a birational
map $f$ is bounded unless the map preserves a fibration. 
If $f$ is holomorphic, it is even bounded independently 
of $f$ (by $dim(M)+h^{1,1}(M)$).

Xie Junyi \cite{X} used this to observe that the orbit closure conjecture is true 
for birational automorphisms of surfaces.

\subsection{Some speculations}

Suppose we want to prove Zhang conjecture for surfaces.
A polarized endomorphism $f:S\to S$, where $S$ is a surface,
admits an ample line bundle $L$ such that $f^*L=L^{\otimes q}$, $q>1$.
The degree of $f$ is then equal to $q^2$. 

We know from holomorphic dynamics that such an endomorphism must have a lot of repelling 
periodic points (they are, in particular, Zariski dense). This is a very non-trivial fact
which follows from the work of Briend-Duval et Dinh-Sibony. The precise condition is that
$deg(f)$ should be greater than the ``first dynamical degree'' which in general is defined
as the limit of $(R(f^n))^{1/n}$, where $R$ denotes the spectral raduis of the action of $f$ on 
$H^{1,1}(S)$ (by inverse image), and here is obviously equal to $q$.

Take one of them and suppose, replacing $f$ by a power if necessary, that this point $x$ is
fixed. Let $\lambda_1, \lambda_2$ denote the eigenvalues of the differential at $x$.
If these are resonant, the resonance should look like $\lambda_1=\lambda_2^a$ where $a\geq 2$ is
an integer (this is since both eigenvalues are greater than one in absolute value). Blowing $x$ up,
we obtain a rational self-map $f':S'\dasharrow S'$, where the fixed point $x$ is replaced by two 
fixed points $x_1, x_2$ on the exceptional divisor. The eigenvalues at $x_1$ are 
$(\lambda_1/\lambda_2, \lambda_2)$, and at $x_2$ they are $(\lambda_1, \lambda_2/\lambda_1)$.
That is, by blowing-up we can eventually produce a non-resonant fixed point with equal
eiganvalues $\lambda_1=\lambda_2=\lambda$ and linearize the map in its neighbourhood. 

If the differential is a Jordan cell, it is an exercize to show that the orbit of almost every
point is Zariski-dense. If the differential is scalar, then in the linearizing coordinates $(y_1, y_2)$ 
we have invariant lines $y_1=\alpha y_2$ for any $\alpha$. Most of those analytic lines would not
``close'' to give algebraic curves, and shall be Zariski-dense in $S$. We must somehow justify that
the analytic lines which do close down to algebraic curves are too ``sparce'' to contain all
algebraic points.

As we have seen, there are examples where infinitely many invariant lines correspond to 
algebraic curves, like the $d$th power map on a rational surface. One also
can produce the 
examples of the same type but with an abelian variety as $S$. The following 
vague question seems
to be natural to ask.

\begin{quest} When there are invariant curves of arbitrarily large degree
through a fixed point of a rational self-map of an algebraic surface, yet
no fibration is preserved, is there always a group structure somewhere behind this
picture?
\end{quest}

Note also that in the case of the $d$th power map, the closure of an 
invariant line
is an algebraic curve if and only if its slope is rational (in suitable
coordinates): $y=x^{\alpha}$ linearizes as $y_2=\alpha y_1$.

An even more vague but also natural question is then as follows:

\begin{quest} Can the invariant algebraic curves through a point be in a
natural correspondence with the algebraic, rather than natural, numbers?
\end{quest}

I believe that any nontrivial information on these two problems would mean a 
major step towards the proof of Zhang's and orbit closure conjectures.

{}

\end{document}